\documentclass{amsart}

\newtheorem{theorem}{Theorem}[section]

\newcommand{\T}{\mathbb{T}}
\theoremstyle{definition}
\newtheorem{definition}[theorem]{Definition}

\theoremstyle{remark}

\numberwithin{equation}{section}

\begin{document}

\title{Accelerated relaxation enhancing flows cause total dissipation}


\author{Keefer Rowan}
\address{Courant Institute of Mathematical Sciences, New York University, New York 10012}
\email{keefer.rowan@cims.nyu.edu}


\date{}

\commby{}

\begin{abstract}
    We show that by ``accelerating'' relaxation enhancing flows, one can construct a flow that is smooth on $[0,1) \times \mathbb{T}^d$ but highly singular at $t=1$ so that for any positive diffusivity, the advection-diffusion equation associated to the accelerated flow totally dissipates solutions, taking arbitrary initial data to the constant function at $t=1$.
\end{abstract}

\maketitle

\section{Introduction}

We consider the evolution of advection-diffusion equations on the torus $\mathbb{T}^d$ with incompressible advecting flow. That is we study solutions to the equation
\[\partial_t \theta - \nu \Delta \theta + u \cdot\nabla \theta =0,\]
where $\nabla \cdot u=0$ and $\nu \geq 0.$ As this equation preserves the mean of $\theta$, we will without loss of generality suppose $\int_{\mathbb{T}^d} \theta(x)\,dx =0$ throughout. The presence of a positive diffusivity $\nu>0$ ensures that $\|\theta\|_{L^2(\T^d)}(t)$ is strictly decreasing in time. A flow is relaxation enhancement if it causes this dissipation of the $L^2$ norm of $\theta$ to happen faster than if no flow were present. To make this precise, let us introduce some definitions from the literature~\cite{constantin_diffusion_2008,feng_dissipation_2019,zelati_relation_2020}. These ideas originate in~\cite{constantin_diffusion_2008}, which originally defined relaxation enhancing. We however use a different---but equivalent---definition of relaxation enhancing.
\begin{definition}
    We denote
    \[L^2_0(\T^d) := \{\theta \in L^2(\T^d) : \int \theta\, dx = 0\}.\]
\end{definition}

\begin{definition}
    For a flow $u \in L_{loc}^\infty([0,\infty) \times \T^d)$ such that $\nabla \cdot u =0$, a diffusivity $\nu > 0$, and times $0\leq s \leq t<\infty$, let $\Phi_{s,t}^{u,\nu} :L^2_0(\T^d) \to L^2_0(\T^d)$ denote the solution operator to the PDE
    \begin{equation}
    \label{eq.advection-diffusion}
    \partial_t \theta - \nu \Delta \theta + u \cdot \nabla \theta =0\end{equation}
    so that $\theta(t) := \Phi_{s,t}^{u,\nu}\theta_0$ solves~\eqref{eq.advection-diffusion} on $[s,\infty) \times \T^d$ with initial data $\theta_0.$
\end{definition}

\begin{definition}
    Define the \textit{dissipation time} of a flow as the maximum time it take for half of the $L^2$ norm of a solution to~\eqref{eq.advection-diffusion} to diffuse, that is
    \[\tau^u(\nu):= \sup_{s\geq 0} \Big(\inf \{t-s : t\geq s, \|\Phi^{u,\nu}_{t,s}\|_{L^2_0 \to L^2_0} \leq \tfrac{1}{2}\}\Big).\]
\end{definition}

\begin{definition}
    Say a flow $u \in L_{loc}^\infty([0,\infty) \times \T^d)$ such that $\nabla \cdot u =0$ is \textit{relaxation enhancing} if
    \[\lim_{\nu \to 0} \nu \tau^u(\nu) =0.\]
\end{definition}

A very wide variety of flows are relaxation enhancing, including all mixing flows, as will be discussed below. Before stating the main result, we need to introduce the notion of accelerating a flow.

\begin{definition}
    For a flow $u \in L^\infty_{loc}([0,\infty) \times \T^d)$ and an increasing diffeomorphism $\sigma: [0,1) \to [0,\infty)$, we define the accelerated flow $u_\sigma \in L^\infty_{loc}([0,1) \times \T^d)$ by
    \[u_\sigma(t,x) := \sigma'(t)u(\sigma(t),x).\]
\end{definition}

This definition is given so that if $\theta$ solves the transport equation associated to $u$, then $\theta(\sigma(t),x)$ solves the transport equation associated to $u_{\sigma}$. We can now state the main result.

\begin{theorem}
    \label{thm.main-result}
    Suppose that $u \in L^\infty_{loc}([0,\infty) \times \T^d)$ such that $\nabla \cdot u =0$ is relaxation enhancing. Then there exists some acceleration of $u$ that is totally diffusive on $[0,1]$ for all diffusivities $\nu>0$, that is there exists some smooth diffeomorphism $\sigma : [0,1) \to [0,\infty)$ such that for all $\nu>0$, $s \in [0,1), \theta_0 \in L^2_0(\T^d),$
    \[\Phi^{u_{\sigma},\nu}_{s,1}\theta_0 =0.\]
\end{theorem}

A natural question is what space $u_\sigma$ belongs to. Let $X$ be some Banach space of functions $\T^d \to \mathbb{R}^d$. For the sake of clarity, we consider the case that the flow $u$ is constant-in-time (a perfectly valid possibility), though this discussion is easily adapted to a variety of assumptions on a time-dependent flow, most simply that $u \in C_b([0,\infty), X).$  

We first note that $u_\sigma \not \in L^1([0,1], X)$, as 
\[\int_0^1 \|u_\sigma\|_X(t)\, dt =\|u\|_X \int_0^1 \sigma'(t)\, dt = \|u\|_X \sigma(1) = \infty.\]
On the other hand, we have that $u_\sigma \in C^\infty([0,1), X).$ Thus we see that $u_\sigma$ is regular away from $t=1$, but so highly singular at $t=1$ so as not to live in any Lebesgue space. On the other hand, $u_\sigma$ lives quite naturally in weighted-in-time spaces, in particular we trivially have
\[\frac{1}{\sigma'} u_\sigma \in L^\infty([0,1],X),\]
where we note that $\lim_{t \to 1} \frac{1}{\sigma'(t)} = 0$. Then to determine which weighted spaces $u_\sigma$ belongs to, we need to get an upper bound for $\sigma'$. The bound on $\sigma'$ in turn depends on how relaxation enhancing the flow $u$ is---that is how fast $\nu \tau^u(\nu) \stackrel{\nu \to 0}{\to} 0$. The next result follows from a general bound on $\sigma'$ and gives the general weighted space that $u_\sigma$ belongs to. Additionally, we give a more concrete bound in the particularly relevant case of exponential mixing.

\begin{theorem}
    \label{thm.bound}
    Let $X$ be some Banach space of functions $\T^d \to \mathbb{R}^d$. Then the $u_\sigma$ in Theorem~\ref{thm.main-result} can be taken so that
    \begin{equation}
        \label{eq.general-bound}
        \Big\|\frac{f^{-1} \Big(\frac{1-t}{4(|\log_2(1-t)| +1)}\Big)}{2 (|\log_2 (1-t)|+1)}u_\sigma(t)\Big\|_{L^\infty([0,1],X)} \leq \|u\|_{L^\infty([0,\infty),X)},
    \end{equation}
    where
    \[f(a) := \sup_{0 \leq b \leq a} b \tau^u(b).\]
    In particular, if $u \in L^\infty([0,\infty), C^{\infty}(\T^d))$ is exponentially mixing, that is if for every $\theta_0 \in L^2_0(\T^d)$, we have for some $K,p>0,$ the estimate
    \begin{equation}
    \label{eq.mixing-est}    
    \|\Phi_{0,t}^{u,0}\theta_0\|_{\dot H^{-1}(\T^d)} \leq K e^{-t^p/K} \|\theta_0\|_{H^1(\T^d)},
    \end{equation}
    then we can take $u_\sigma$ in Theorem~\ref{thm.main-result} so that
    \begin{equation}
    \label{eq.exp-mixer-est}    
    \frac{1-t}{|\log_2(1-t)|^{2+2/p} + 1}u_\sigma(t) \in L^\infty([0,1],C^\infty(\T^d)).
    \end{equation}
\end{theorem}

\section{Discussion}

Recently, there has been substantial interest in the phenomenon of anomalous dissipation for the passive scalar advection-diffusion equation~\cite{drivas_anomalous_2022,colombo_anomalous_2023,armstrong_anomalous_2023,elgindi_norm_2023,burczak_anomalous_2023}. In these works, an incompressible flow $u \in L^\infty([0,1] \times \T^d)$ is constructed\footnote{In each work, $u$ belongs to much stronger space than just $L^\infty_{t,x}$, but for our sake this is the relevant fact.} such that for some initial data $\theta_0 \in L^2(\T^d)$,\footnote{In most of these works, a substantially stronger result is shown than anomalous dissipation for just \textit{some} initial data. In particular, in~\cite{drivas_anomalous_2022}, anomalous dissipation is shown for all initial data that is sufficiently close to eigenfunctions of the Laplacian. In~\cite{armstrong_anomalous_2023}, anomalous dissipation is shown for all $\theta_0 \in H^1(\T^d)$. In~\cite{elgindi_norm_2023}, anomalous dissipation is shown for all $\theta_0 \in W^{1,\infty}(\T^d) \cap H^{1+s}(\T^d)$ for some $s>2/5.$ In~\cite{burczak_anomalous_2023}, the construction of~\cite{armstrong_anomalous_2023} is modified so that the flow $u$ solves the Euler equation.}
\[\liminf_{\nu \to 0} \|\Phi^{u,\nu}_{0,1} \theta_0\|_{L^2(\T^d)} < \|\theta_0\|_{L^2(\T^d)}.\]
In Theorem~\ref{thm.main-result}, we show total dissipation, $\Phi_{0,1}^{u,\nu} \theta_0 = 0$, uniform in diffusivity. As such, this is a special case of anomalous dissipation, though the flow $u$ we construct is much less physical than these other examples of anomalous dissipation, as it is highly singular at the final time.

In~\cite{hofmanova_anomalous_2023}, total dissipation for any diffusivity $\nu \in (0,1)$ is shown for a scalar advected by a solution to a randomly forced Navier-Stokes equation. The flow used in~\cite{hofmanova_anomalous_2023} is also highly singular at the final time, failing to belong to any $L^p$ space on the full interval $(0,1)$. Theorem~\ref{thm.main-result} shows that the flows of~\cite{hofmanova_anomalous_2023} are part of a broad class of highly singular flows that cause total dissipation. It is worth noting that total dissipation is impossible at any finite diffusivity $\nu >0$ if the flow $u \in L^\infty_{t,x}$, by the unique continuation result of \cite{poon_unique_1996}.\footnote{Note however that this doesn't prevent one from having asymptotically total dissipation in the limit as $\nu \to 0$, i.e.\ that $\lim_{\nu \to 0} \Phi^{u,\nu}_{0,1} \theta_0 = 0$, for a flow $u \in L^{\infty}_{t,x}$.}

Lastly, we note that~\cite{feng_dissipation_2019,zelati_relation_2020} show that any mixing flow is relaxation enhancing.\footnote{For precise definitions of the relevant sense of mixing, see ~\cite{feng_dissipation_2019,zelati_relation_2020}.} In particular,~\cite{feng_dissipation_2019,zelati_relation_2020} show quantitative relations between mixing rates and dissipation times. This implies the broad class of mixing flows are relaxation enhancing and further gives a quantitative bound on the rate of relaxation enhancement (the rate that $\nu \tau^u(\nu) \to 0$) in terms of the rate of mixing. Thus through Theorem~\ref{thm.bound}, quantitatively mixing flows allow us to construct $u_\sigma$ which belong to specific weighted $L^\infty$ spaces. As given in Theorem~\ref{thm.bound}, we compute the weighted space $u_\sigma$ belongs to when $u$ is exponentially mixing, but this computation can be repeated for any given rate of mixing, using~\cite{feng_dissipation_2019} to convert the mixing rate to an enhanced dissipation rate and then using Theorem~\ref{thm.bound} to determine the appropriate weighted space. The existence of an exponential mixer $u \in L^\infty([0,\infty),C^\infty(\T^d))$ with $p=1$ is given by~\cite{blumenthal_exponential_2023}, so there are flows $u$ for which we can apply the estimate~\eqref{eq.exp-mixer-est}.

\section{Proofs}

We now provide the straightforward and short proofs of Theorems~\ref{thm.main-result} and~\ref{thm.bound}.

\begin{proof}[Proofs of Theorems~\ref{thm.main-result} and~\ref{thm.bound}]
    We define the diffeomorphism $\sigma$ as a regularization of a piecewise linear flow. Let us first specify the piecewise linear flow. Fix some strictly increasing sequence of time $T_j$ such that $T_0 =0$; the remaining $T_j$ will be specified later. Then let
    \[\tilde \sigma(1-2^{-j}) = T_j,\]
    with $\tilde \sigma(t)$ taken to be linear on $[1-2^{-j}, 1-2^{-(j+1)}].$ We then define $\sigma$ as a strictly increasing regularization of $\tilde \sigma$ so that $\sigma \in C^\infty(0,1)$ and so that $\sigma' \leq 2\tilde \sigma'$ and for any $j \in \mathbb{N}$ and and any $t \in [1 - 2^{-(j-1)}, 1-2^{-(j-1)} + 2^{-(j+1)}]$
    \[\sigma(t) = \tilde \sigma(t).\]
    Then $\sigma : [0,1) \to [0,\infty)$ is a smooth diffeomorphism. We now fix $\nu>0$. Our goal now is to choose $T_j$ indepedently of $\nu$ so that for $j$ sufficiently large
    \[\|\Phi^{u_\sigma,\nu}_{1 - 2^{-(j-1)}, 1-2^{-(j-1)} + 2^{-(j+1)}}\|_{L^2_0(\T^d) \to L^2_0(\T^d)} \leq \tfrac{1}{2}.\]
    Note that this clearly implies the total dissipation of Theorem~\ref{thm.main-result}.
    
    By a simple change of variables
    \[\Phi^{u_\sigma,\nu}_{1 - 2^{-(j-1)}, 1-2^{-(j-1)} + 2^{-(j+1)}} =\Phi^{u, \frac{2^{-j} \nu}{T_{j} -T_{j-1}}}_{T_{j-1},\frac{T_{j-1}+T_{j}}{2}},\]
    thus by the definition of $\tau^u$, it suffices to choose $T_j$ so that for $j$ sufficiently large
    \[\tau^u\Big(\frac{2^{-j} \nu}{T_{j} -T_{j-1}}\Big) \leq \frac{T_{j}-T_{j-1}}{2}.\]
    This in turn is implied by 
    \[2^{-(j+1)} \nu \geq f\Big(\frac{2^{-j} \nu}{T_{j} -T_{j-1}}\Big),\]
    where $f$ is defined as in Theorem~\ref{thm.bound}. Note that $f$ is increasing and for any $\nu>0$, we have that eventually $j^{-1} \leq \nu \leq j$, so it suffices to choose $T_j$ so that
    \begin{equation}
        \label{eq.T_j-def}
    \frac{2^{-(j+1)}}{j} = f\Big(\frac{2^{-j} j}{T_{j} -T_{j-1}}\Big).
    \end{equation}
    Thus we take $T_j$ so that
    \[T_j - T_{j-1} := \frac{2^{-j} j}{f^{-1}\Big(\frac{2^{-(j+1)}}{j}\Big)}.\]
    Choosing $T_j$ in this way concludes the construction of $\sigma$ and by the arguments above concludes the proof of Theorem~\ref{thm.main-result}. 
    
    What remains is to prove the bounds of Theorem~\ref{thm.bound}. Note that on the interval $t \in (1-2^{-(j-1)}, 1-2^{-j})$, we have
    \[\sigma'(t) \leq 2 \tilde\sigma'(t) = 2^{j+1} (T_j - T_{j-1}) = \frac{2j}{f^{-1}\Big(\frac{2^{-(j+1)}}{j}\Big)}.\]
    Note then that on this interval
    \[2^{-j} \leq 1 - t  \leq 2^{-(j-1)},\]
    thus using that $f^{-1}$ is increasing
    \[\sigma'(t) \leq \frac{2 (|\log_2 (1-t)|+1)}{f^{-1} \Big(\frac{1-t}{4(|\log_2(1-t)| +1)}\Big)}.\]
    Plugging this into the definition of $u_\sigma$, we get~\eqref{eq.general-bound}.

    In order to conclude, we just need to specialize~\eqref{eq.general-bound} to the case of exponential mixing. Suppose now that $u$ satisfies the mixing estimate~\eqref{eq.mixing-est}. From~\cite{zelati_relation_2020},~\eqref{eq.mixing-est} implies that there exists some $C(K,p)<\infty$ such that for $\nu \leq C^{-1}$
    \[\tau^u(\nu) \leq C |\log \nu|^{2/p}.\]
    Thus for $a \leq C^{-1},$
    \[f(a) \leq Ca |\log a|^{2/p}.\]
    It's somewhat unwieldy to compute a good bound directly from~\eqref{eq.general-bound}, so let us instead return to~\eqref{eq.T_j-def} to give for $j$ sufficiently large,
    \begin{align*} \frac{2^{-(j+1)}}{j} &= f\Big(\frac{2^{-j} j}{T_{j} -T_{j-1}}\Big) 
    \\&\leq C \frac{2^{-j} j}{T_{j} -T_{j-1}} \bigg| \log \frac{2^{-j} j}{T_{j} -T_{j-1}} \bigg|^{2/p} 
    \\&\leq C\frac{2^{-j} j}{T_{j} -T_{j-1}}\Big(j^{2/p} + |\log (T_j-T_{j-1})|^{2/p}\Big).
    \end{align*}
    Thus
    \[T_j - T_{j-1} \leq C j^2 \Big(j^{2/p} + |\log (T_j-T_{j-1})|^{2/p}\Big).\]
    One can verify that this inequality implies that
    \[T_j - T_{j-1} \leq C j^{2+2/p}.\]
    Then, as above, we have that on the interval $t \in (1-2^{-(j-1)}, 1-2^{-j}),$
    \[\sigma'(t) \leq 2 \tilde \sigma'(t) = 2^{j+1} (T_j - T_{j-1}) \leq C 2^{j+1} j^{2+2/p} \leq C\frac{|\log_2(1-t)|^{2+2/p} + 1}{1-t}.\]
    This estimate then gives the result.
\end{proof}

\textit{Acknowledgements.} I would like to thank Scott Armstrong and Vlad Vicol for stimulating discussion. The author was partially supported by NSF grants DMS-1954357 and DMS-2000200 as well as a Simons Foundation grant.

\bibliographystyle{amsplain}
\bibliography{references}

\end{document}